\documentclass[twocolumn,aps,floatfix,superscriptaddress,longbibliography]{revtex4-1}
\pdfoutput=1 
\usepackage{amsmath,amssymb,eucal,graphicx,float,epstopdf,xparse}
\usepackage{epsfig,subfigure}
\usepackage[utf8]{inputenc}

\usepackage[colorlinks=true, urlcolor=blue, anchorcolor=blue, citecolor=blue,filecolor=blue,linkcolor=blue,menucolor=blue]{hyperref}

\begin{document}
\title{Simple evolving random graphs}

\author{P. L. Krapivsky}
\affiliation{Department of Physics, Boston University, Boston, Massachusetts 02215, USA}
\affiliation{Santa Fe Institute, Santa Fe, New Mexico 87501, USA}

\begin{abstract} 
We study the evolution of graphs densifying by adding edges: Two vertices are chosen randomly, and an edge is (i) established if each vertex belongs to a tree; (ii) established with probability $p$ if only one vertex belongs to a tree; (iii) an attempt fails if both vertices belong to unicyclic components. Emerging simple random graphs contain only trees and unicycles. In the thermodynamic limit of an infinite number of vertices, the fraction of vertices in unicycles undergoes a phase transition resembling a percolation transition in classical random graphs. In classical random graphs, a complex giant component born at the transition eventually engulfs all finite components and densifies forever. The evolution of simple random graphs freezes when trees disappear. We quantify simple random graphs in the supercritical phase and the properties of the frozen state. 
\end{abstract}

\maketitle

\section{Introduction}

Classical evolving random graphs \cite{ER60} are built by drawing edges between pairs of randomly chosen vertices. In the most interesting sparse regime, when the number of edges is comparable with the number of vertices, components (i.e., maximal connected sub-graphs) are predominantly trees. Unicyclic components (shortly unicycles) occasionally appear, yet the number of unicycles remains finite even in  the thermodynamic limit when the size of the graph (the number of vertices $N$) diverges. The only exception is the percolation point where the number of unicycles is infinite; when $N\gg 1$, the transition occurs in a scaling window where the number of unicycles scales logarithmically with $N$. The giant component arising in the supercritical  phase is neither a tree nor a unicycle \cite{Janson93}. 

Classical random graphs constitute a developed and growing field of research with applications to computer science, mathematics, and natural sciences \cite{Bollobas,Hofstad,Frieze,KRB,chatterjee16}.  The probabilistic treatment of classical random graphs is popular in mathematical literature \cite{Bollobas,Hofstad,Frieze}. A dynamical treatment was already used by Erd\H{os} and R\'{e}nyi \cite{ER60}. Earlier studies of gelation by Flory and Stockmayer employed other techniques and explored somewhat different random graph models \cite{Flory41,Stockmayer43,Flory53}, see \cite{Ziff80,Leyvraz-rev} for comparison of the Flory and Stockmayer approaches. We employ a dynamical treatment, a version that we call a kinetic theory approach \cite{KRB,Lushnikov05}. Originally, such an approach has been applied to finding the distribution of trees \cite{Mcleod62a,Mcleod62b}. The distribution of the average number of unicycles is also amenable to the kinetic theory treatment \cite{BK-04,BK-05}. 

We distinguish the complexity of components by their Euler characteristics.  For a graph with  $V$ vertices and $E$ edges, the Euler characteristic is 
\begin{equation}
\label{Euler}
\chi=V-E
\end{equation}
The Euler characteristic is a topological invariant. The identities $\chi(\text{tree})=1$ and $\chi(\text{unicycle})=0$ illustrate the topological nature of the Euler characteristic. Trees and unicycles are known as simple components; components with $\chi<0$ are complex \cite{Janson93a}. Simple graphs, by definition, are disjoined unions of trees and unicycles. 

Bicyclic components are born when unicycles merge or after adding an edge in a unicycle. In a classical random graph process, such events are rare. The giant component appearing in the post-percolation phase is complex. The probability that, throughout the evolution, there never is more than a single complex component is $5\pi/18$. Up to the percolation point, the evolving graph consists entirely of trees, unicycles, and bicycles with probability $\sqrt{2/3}\cosh\big(\sqrt{5/18}\big)=0.932\,548\ldots$. See \cite{Janson93a,Stepanov88,Knuth89,Janson93} for these and other subtle behaviors. Overall, a few complex components arise in the classical random graph process. This property suggests exploring graph evolutions in which the formation of complex components is strictly forbidden---attempts to draw edges between unicycles and inside the same unicycle are rejected, so components are either trees or unicycles. 

More precisely we examine a one-parametric class of simple random graph (SRG) processes defined as follows. A disjoint graph with $N$ vertices begins to densify by adding edges according to the following procedure:
\begin{enumerate}
\item An attempt to draw an edge between vertices belonging to trees is always successful, as in the classical random graphs.
\item An attempt to draw an edge between a vertex from a tree and a vertex from a unicycle is successful with gluing probability $p$. 
\item An attempt to draw an edge between vertices from unicycles is rejected.
\end{enumerate}
Since unicycles never merge, it is natural to postulate that drawing an edge between vertices from a tree and a unicycle occurs with smaller rate than between vertices from trees. This suggests to consider the range $p\in [0,1]$ and interpret $p$ as a gluing probability. Mathematically, an SRG process with any $p\geq 0$ is well-defined if we interpret $p$ as a rate, with the rate of drawing an edge between vertices belonging to trees set to unity. 

We defined SRG processes by constraining the random graph process to be free of complex components. There are other constrained random graph processes in which the formation of some subgraphs is forbidden \cite{Bollobas00,Bohman10,Frieze}. In these processes, the evolution eventually stops, and the graph freezes. In the triangle-free process \cite{EW95,Kim95,Bohman09,Morris}, every new edge is randomly chosen from those that are not present and that do not create a triangle. The triangle-free process provides a lower bound for the Ramsey number $R(3,k)$, see \cite{Kim95,Bohman09,Morris}. The triangle-free process ends with a dense frozen graph \cite{Bohman15,Razborov22}; sparse graphs emerge at the end of the SRG processes. The sparsity is a simplifying feature; constrained random graph processes creating dense jammed graphs tend to be more challenging for analyses.

Random graphs with various constraints have been mostly studied in static frameworks \cite{Razborov08,Radin13,Radin14,Razborov17}. We employ an evolution framework.

We now outline a few features of SRGs derived below. The fraction of mass $s$ in unicycles in the supercritical phase, $t>1$, is implicitly determined by 
\begin{equation}
\label{T:sol}
t = 1+\frac{1}{p\,s^\frac{1}{p}}\int_0^s dz\,\frac{z^\frac{1}{p}}{1-z}
\end{equation}
Equation \eqref{T:sol} admits a trivial solution $s=0$ and a non-trivial solution giving the correct value of $s(t)$ in the supercritical  phase. Extracting asymptotic behaviors of $s(t)$ from the implicit solution \eqref{T:sol} yields
\begin{equation}
\label{s:extreme}
s=
\begin{cases}
0                                                                                        & t<1\\
(1+p)(t-1)+\mathcal{O}\big[(t-1)^2\big]                              & t\downarrow 1\\
1 - \mathcal{E} + \mathcal{O}\big[t \mathcal{E}^2\big]      & t\uparrow \infty
\end{cases}
\end{equation}
with
\begin{equation}
\label{E:def}
\mathcal{E}(t) = e^{-pt + p -\gamma-\psi(1+1/p)}
\end{equation}
Here $\gamma=0.577215\ldots$ is the Euler-Masceroni constant, $\psi(z)=\Gamma'(z)/\Gamma(z)$ the digamma function, and $\Gamma(z)$ the gamma function.

In the supercritical  phase, SRGs are significantly different from classical random graphs. In Sec.~\ref{sec:simple}, we show that SRGs undergo a continuous phase transition, determine the fraction of mass in unicycles, and derive the distribution of trees. 

The classical random graph becomes connected at time $t_\text{cond}\simeq \ln N$ and the densification continues ad infinitum.  The SRG process, in contrast,  reaches a jammed state once trees disappear. Using heuristic arguments we estimate the jamming time 
\begin{equation}
\label{jamming}
t_\text{jam} \simeq p^{-1} \ln N
\end{equation}
The average number of unicycles $U(t)$ increases until the system reaches a jammed state where
\begin{equation}
\label{U-jam}
U_\text{jam} \simeq \frac{1+p}{6p}\,\ln N
\end{equation}

Equations \eqref{T:sol}--\eqref{U-jam} do not apply when $p=0$, i.e., to the model with frozen unicycles. Some results simplify in this extreme case, e.g., $s=1-1/t$ in the supercritical phase, $t>1$. The logarithmic scaling laws \eqref{jamming}--\eqref{U-jam} valid for models with $p>0$ are replaced by algebraically scaling behaviors (Sec.~\ref{subsec:extreme-0}). 

In Sec.~\ref{sec:classical}, we outline several properties of classical random graphs needed for the analysis of SRGs. We also probe less known features like the complexity of the giant component. In Secs.~\ref{sec:simple}--\ref{sec:extreme}, we describe the phase transition in the SRGs, determine the densities of trees and the average number of unicycles in the supercritical phase, and probe the properties of the final jammed state. In Sec~\ref{sec:fluct}, we step away from average characteristics and emphasize that fluctuations that could be particularly important for the SRG process with frozen unicycles ($p=0$). We finish with concluding remarks (Sec.~\ref{sec:remarks}). 

\section{Classical Random Graphs}
\label{sec:classical}

It is customary to begin with a disjoint graph, i.e., a collection of $N$ isolated vertices. Edges are drawn between any two vertices with rate $1/(2N)$. With this convention, the percolation transition happens at time $t_c=1$ If $N\gg 1$, almost all components are trees. The number $\mathfrak{T}_k$ of trees of size $k$ is a self-averaging random quantity: $\mathfrak{T}_k = N c_k + \sqrt{N} \xi_k$ with random $\xi_k=O(1)$. Therefore, fluctuations are relatively small, so the deterministic densities $c_k$ provide the chief insight. In the $N\to\infty$ limit, the densities evolve according to an infinite set of coupled non-linear ordinary differential equations (ODEs)
\begin{equation}
\label{prod}
\frac{dc_k}{dt} = \frac{1}{2}\sum_{i+j=k} ij\;c_i c_j - k c_k
\end{equation}
from which \cite{KRB,Leyvraz-rev}
\begin{equation}
\label{prod-densities}
c_k(t)= \frac{k^{k-2}}{k!}\, t^{k-1}\, e^{-kt}
\end{equation}
The smoothness of densities throughout the entire evolution, $0<t<\infty$, hides the emergence of the giant component when $t>t_c=1$. Mass conservation, $\sum kc_k(t)=1$, is manifestly obeyed in the subcritical (pre-percolation) phase. In the supercritical  phase, $\sum kc_k(t)=1-g(t)$, where $g(t)$ is the fraction of vertices belonging to the giant component, shortly the mass of the giant component. The mass $g(t)$ is implicitly determined by \cite{KRB,Leyvraz-rev}
\begin{equation}
\label{prod-giant}
g=1-e^{-gt}
\end{equation}
For instance, the giant component comprises half mass of the system, $g=1/2$, at time $t=2\ln 2$. 

The total cluster density $c(t)=\sum_{k\geq 1} c_k(t)$ reads
\begin{equation}
\label{prod-total}
c(t) =
\begin{cases}
1-\frac{t}{2}                                                          & t\leq 1\\
1-\frac{t}{2} +\frac{1}{2}\int_1^t d\tau\,g^2(\tau)   & t>1
\end{cases}
\end{equation}
The integral can be expressed as a function of $g$:
\begin{eqnarray}
\label{integral}
\int_1^t g^2(\tau)\, d\tau &=& \int_0^g h^2\, d[-h^{-1}\ln(1-h)] \nonumber \\
&=& \int_0^g \left[\ln(1-h) + \frac{h}{1-h}\right]dh \nonumber\\
&=& -2g -(2-g)\ln(1-g)
\end{eqnarray}
Here $h=g(\tau)$ and we use $\tau=-h^{-1}\ln(1-h)$ following from \eqref{prod-giant}. Thus in the supercritical  phase ($t>1$)
\begin{equation}
\label{total-post}
c(t) = 1-\frac{t}{2} -g-\left(1-\frac{g}{2}\right)\ln(1-g)    
\end{equation}
Using this formula we extract more explicit results for the cluster density $c(t)$ just beyond the critical point and in the long time limit:
\begin{equation}
\label{total-more}
c =
\begin{cases}
1-t/2+2(t-1)^3/3+\ldots      & t\downarrow 1\\
e^{-t}+(t/2)\,e^{-2t}+\ldots & t\uparrow \infty
\end{cases} 
\end{equation}

In the following we shall also need the second moment $M_2 = \sum_{k\geq 1} k^2 c_k$. It can be expressed \cite{KRB} as follows
\begin{equation}
\label{prod-M2}
M_2(t)=
\begin{cases}
(1-t)^{-1}                  &  t<1\\
\frac{1-g}{1-t(1-g)}   &   t>1
\end{cases} 
\end{equation}
In the supercritical (post-percolation) phase, the second moment accounts only for finite components. (In the leading order, the second moment is $Ng^2$ if we include the contribution of the giant component.) The second moment rapidly decays with time as the giant component quickly engulfs finite components.

The giant component has more edges than vertices. The ratio of the number of edges to the number of vertices in the giant component is \cite{BK-04}
\begin{equation}
\label{prod-ratio}
\frac{\#(\text{edges})}{\#(\text{vertices})} = 1 + \frac{1}{2g(t)}\int_1^t d\tau\,g^2(\tau)
\end{equation}
where the integral accounts that the number of edges within the giant component increases with rate $g^2/2$. The integral on the right-hand side of  \eqref{prod-ratio} is given by \eqref{integral}, so Eq.~\eqref{prod-ratio} simplifies to
\begin{equation}
\label{ratio}
\frac{\#(\text{edges})}{\#(\text{vertices})} = -\frac{2-g}{2g}\,\ln(1-g)
\end{equation}
By definition, the number of  vertices in the giant component is $\#(\text{vertices})=gN$, and from \eqref{Euler} and \eqref{ratio} we find the Euler characteristic of the giant component:
\begin{eqnarray}
\chi_\text{giant} &=& N\left[g+\frac{2-g}{2}\,\ln(1-g)\right]\nonumber  \\
&=& -N\sum_{m\geq 3} \frac{m-2}{2m(m-1)}\,g^m
\end{eqnarray}
The giant component is complex as its Euler characteristic is negative. Further, the Euler characteristic of the giant component is giant, viz. extensive in system size. 

\begin{figure}[t]
\includegraphics[width=0.44\textwidth]{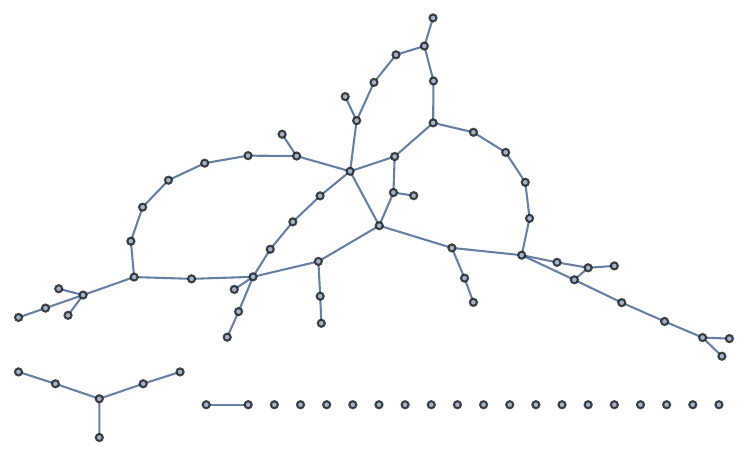}
\caption{An example of a random graph with $N=80$ vertices and $E=65$ edges. The average number of edges in an evolving random graph is $Nt/2$. Hence, it is natural to interpret the above graph as a possible outcome of the random graph process at time $t=2E/N=13/8$, and compare with theoretical predictions in the supercritical phase at $t=13/8$.} 
\label{fig:CRG-post}
\end{figure}

As an illustration, in Fig.~\ref{fig:CRG-post}, we show a random graph with $N=80$ vertices and $E=65$ edges. Such a graph is generated, on average, at time $t=13/8$. So, the graph is in the supercritical phase, and indeed, it has a complex component, $\chi=-5$, that is significantly larger than other components. The fraction of mass in the giant component of the graph in Fig.~\ref{fig:CRG-post} is $\frac{54}{80}=0.675$, a little larger than the expected $g\approx 0.655$ predicted by \eqref{prod-giant} and realized when $N\to\infty$. The average total number of isolated vertices is $N_1=N  e^{-t}\approx 15.753$; the graph shown in Fig.~\ref{fig:CRG-post} has $N_1=18$ isolated vertices. From \eqref{Ut-sol-post} we get $U\approx 0.41$; the graph in Fig.~\ref{fig:CRG-post} has no finite unicycles. The ratio of the number of edges to the number of vertices in the giant component is $\frac{59}{54}=1.0925925\ldots$ for the graph shown in Fig.~\ref{fig:CRG-post}; the predicted ratio \eqref{ratio} realized in the $N\to\infty$ limit is $1.0927078\ldots$.

Unicycles in classical random graphs have been investigated using probabilistic and combinatorial techniques \cite{Stepanov88,Knuth89,Janson93a,Janson93,Janson03}, kinetic theory approach \cite{BK-04,BK-05}, and other methods, some of which are more suitable to random graph models different from classical random graphs \cite{Marinari04,Burda04a,Burda04b,Ginestra05}. The number $\mathfrak{U}_k$ of unicycles of size $k$ is a non-self-averaging random quantity, i.e., fluctuations are comparable with the average. Still, the average $U_k=\langle \mathfrak{U}_k\rangle$ sheds light on unicycles. The kinetic theory approach \cite{BK-04,BK-05} leads to the following infinite set of ODEs 
\begin{equation}
\label{Ukt}
\frac{dU_k}{dt}=\frac{1}{2}\,k^2c_k+\sum_{i+j=k}iU_i\,jc_j-k\,U_k
\end{equation}
These equations are linear and, in this respect, simpler than Eqs.~\eqref{prod} for the densities of trees. An extra `source' term, the first term on the right-hand side, accounting for turning a tree into a unicycle by creating an edge inside the tree, is a complication. 

Summing Eqs.~\eqref{Ukt} we deduce an equation
\begin{equation}
\label{U-M2}
\frac{dU}{dt}=\frac{1}{2}M_2
\end{equation}
for the average total number of unicycles $U=\sum U_k$ in the subcritical  phase. Using $M_2=1/(1-t)$ we find
\begin{subequations}
\begin{equation}
\label{Ut-sol-pre}
U(t)=\frac{1}{2}\,\ln \frac{1}{1-t}
\end{equation}
in the subcritical  phase. Similarly, the average number of finite unicycles is \cite{BK-04,BK-05}
\begin{equation}
\label{Ut-sol-post}
U(t)=\frac{1}{2}\,\ln \frac{1}{1-(1-g)t}
\end{equation}
\end{subequations}
in the supercritical  phase.

For a finite system, the critical point broadens into the scaling window, $|1-t|\sim N^{-1/3}$.  Using this and either \eqref{Ut-sol-pre} or \eqref{Ut-sol-post}  we find that the total number of unicycles remains finite but scales logarithmically with the total number of vertices in the critical regime
\begin{equation}
\label{U-gel}
U(1)=\tfrac{1}{6}\ln N
\end{equation}

Solving \eqref{Ukt} subject to $U_k(0)=0$ recurrently for small $k$ one guesses 
\begin{equation}
\label{Ukt-ansatz}
U_k(t)= B_k\, t^k\,e^{-kt}
\end{equation}
Substituting \eqref{Ukt-ansatz} into \eqref{Ukt} one fixes the amplitudes:
\begin{equation}
\label{Ukt-sol}
U_k(t)=\frac{1}{2}\,t^k\,e^{-kt}\,\sum_{n=0}^{k-1}\frac{k^{n-1}}{n!}
\end{equation}
At the critical point, the average number of unicylces decays algebraically with size of the cycle,  
\begin{equation}
\label{Uk-gel}
U_k(1) =  \frac{1}{2}\,e^{-k}\,\sum_{n=0}^{k-1}\frac{k^{n-1}}{n!}\simeq \frac{1}{4k}
\end{equation}
for $k\gg 1$. For finite graphs, the largest components have size $O(N^{2/3})$ at the critical point \cite{Janson93,Bollobas,Hofstad,Frieze,KRB}. Summing \eqref{Uk-gel} up to $k\sim N^{2/3}$ one obtains \eqref{U-gel}.

In addition to the distribution $U_k(t)$ of unicycles, one can determine the distribution of cycles \cite{Janson93} 
\begin{equation}
\label{C-ell-sol}
C_\ell(t) = \frac{t^\ell}{2\ell}
\end{equation}
of length $\ell$. At the critical point, $C_\ell(1) = \frac{1}{2\ell}$ is strictly algebraic, with a crossover at $\ell\sim N^{1/3}$ for finite graphs. 
The joint distribution of unicycles of size $k$ with cycle of length $\ell\leq k$ is also known \cite{BK-05}
\begin{equation}
\label{Uk-ell-sol}
U_{k,\ell}(t)=\frac{1}{2}\,\frac{k^{k-\ell-1}}{(k-\ell)!}\,t^k\,e^{-kt}
\end{equation}
Using this result one can re-derive \eqref{Ukt-sol} and \eqref{C-ell-sol}. Another consequence of \eqref{Uk-ell-sol} is the formula $R_\ell = U_{\ell,\ell}=\frac{1}{2\ell} t^\ell e^{-\ell t}$ for the average number of rings of length $\ell$. 

We now turn to the SRG evolution processes in which, by definition, complex components cannot be created.

\section{Simple Random Graphs}
\label{sec:simple} 

In the SRG evolution processes, the merging of unicycles is forbidden. The merging of two trees, say a tree of size $i$ and a tree of size $j$, proceeds at the same rate
\begin{subequations}
\begin{equation}
\label{TTT}
[T_i] \oplus [T_j] \to [T_{i+j}] \qquad \text{rate}\quad ij
\end{equation}
as for the classical random graph process. The merging of a tree and a unicycle proceeds via
\begin{equation}
\label{TUU}
[T_i] \oplus [U_j] \to [U_{i+j}] \qquad \text{rate}\quad pij
\end{equation}
where $p\geq 0$. Since unicycles cannot merge, it is natural to assume that the reaction channel \eqref{TUU} is slower than \eqref{TTT}, that is, $0\leq p\leq 1$. This agrees with interpretation that an attempt to draw an edge between a vertex from a tree and a vertex from a unicycle is successful with probability $p$. Finally, drawing an edge inside a tree turns it into a unicycle
\begin{equation}
\label{TU}
[T_k] \to [U_k] \qquad \text{rate}\quad k^2/(2 N)
\end{equation}
\end{subequations}
Two unicycles never merge and drawing an edge inside a unicycle is also forbidden---such processes would have generated complex components. Thus the total number of unicycles increases via the reaction channel \eqref{TU}.

\subsection{Trees and the phase transition}

The influence of the reaction channel \eqref{TUU} on the evolution of trees is asymptotically negligible in the subcritical  phase since the number of unicycles is $N$ times smaller than the number of trees. Therefore, Eqs.~\eqref{prod} are applicable as in the classical random graph, and the densities of trees are given by the same Eq.~\eqref{prod-densities} when $t\leq 1$. Below, we show that the influence of unicycles on the evolution of trees is also negligible in the supercritical  phase. Only in the proximity of the jammed state unicycles affect the evolution of trees. 

In classical random graphs, the components have a size up to $O(N^{2/3})$ at the critical point. The maximal merging rate is $N^{2/3}\times N^{2/3}/N=N^{1/3}$ explaining why the width of the scaling window is $|1-t|\sim N^{-1/3}$. When $t-1\gg N^{-1/3}$, the most massive components merge, and a single giant component emerges. This giant component progressively engulfs finite components. Eventually, only the giant component remains.

In SRGs, trees of size $O(N^{2/3})$ merge with the rate $O(N^{1/3})$, and any such tree turns into a unicycle with the same rate. Thus, trees of size $O(N^{2/3})$ become unicycles in the scaling window. 

The nature of the phase transition is captured by the behavior of $s(t)=N^{-1}\sum_{k\geq 1} k \mathfrak{U}_k(t)$, the fraction of vertices belonging to unicycles, shortly the mass of unicycles. This quantity plays the role of an order parameter. For finite systems $s(t)\sim N^{-1}$ in the subcritical  phase and $s\sim N^{-1/3}$ in the critical regime [see  Sec.~\ref{subsec:uni} for details]. In the supercritical phase, $s(t)$ is finite and growing with time. 
 
The densities of trees satisfy
\begin{equation}
\label{post-long}
\frac{dc_k}{dt} = \frac{1}{2}\sum_{i+j=k} ij\;c_i c_j - k c_k\left(\sum_{\ell\geq 1} \ell c_\ell+p s\right)
\end{equation}
The second sum on the right-hand side of \eqref{post-long} goes over all trees, so 
\begin{equation}
\label{def:s}
\sum_{\ell\geq 1} \ell c_\ell=1-s
\end{equation}
by the definition of $s$. Thus Eq.~\eqref{post-long} simplifies to
\begin{equation}
\label{post}
\frac{dc_k}{dt} = \frac{1}{2}\sum_{i+j=k} ij\;c_i c_j - k c_k\left(1-qs\right)
\end{equation}
with $q=1-p$. In the subcritical  phase, $s=0$ in the leading order, so \eqref{post} and \eqref{prod} coincide, and the solution \eqref{prod-densities} remains valid. Finding $s(t)>0$ in the supercritical  phase is part of the solution. 

We solve Eqs.~\eqref{post} by employing the same approach \cite{KRB} as for classical random graphs. The emerging solution depends on yet-unknown $s$, but by inserting this formal solution into relation \eqref{def:s}, we will fix $s$. The approach used in \cite{KRB} relies on the generating function technique \cite{Flajolet}. It is convenient to introduce an exponential generating function based on the sequence $k c_k$: 
\begin{equation}
\label{C:def}
\mathcal{C}(y,t) = \sum_{k\geq 1} k\,c_k(t)\, e^{yk}
\end{equation}
This allows us to recast an infinite set of ODEs, Eqs.~\eqref{post}, into a single partial differential equation 
\begin{equation}
\label{C:eq}
\partial_t \mathcal{C} = (\mathcal{C}-1+q s)\partial_y \mathcal{C}
\end{equation}
Let us consider $y = y(\mathcal{C},t)$ instead of $\mathcal{C}=\mathcal{C}(y,t)$. This transformation recasts \eqref{C:eq} into
\begin{equation}
\label{y:eq}
\partial_t y +\mathcal{C}-1+qs = 0
\end{equation}
which is integrated to give
\begin{equation}
\label{y:sol-f}
y +(\mathcal{C}-1)t+q\int_0^t d\tau\,s(\tau) = f(\mathcal{C})
\end{equation}
Initially $c_k(0)=\delta_{k,1}$, and therefore $\mathcal{C}(y,0) = e^{y}$. This fixes $f(\mathcal{C})=\ln \mathcal{C}$, and \eqref{y:sol-f} leads to 
\begin{equation}
\label{y:sol}
\mathcal{C}\,e^{-\mathcal{C} t} = e^{y-t+q\int_0^t d\tau\,s(\tau)}
\end{equation}

We re-write Eq.~\eqref{y:sol} as $Z e^{-Z} = Y$ with $Z = \mathcal{C} t$ and $Y = t e^{y-t+q\int_0^t d\tau\,s(\tau)}$. 
The power series $Y=Y(Z)$  directly follow from $Z e^{-Z} = Y$. The inverse power series $Z=Z(Y)$ can be deduced using the Lagrange inversion formula. One gets \cite{KRB,Flajolet}
\begin{equation}
\label{Lagrange}
Z = \sum_{k\geq 1} \frac{k^{k-1}}{k!}\,Y^k
\end{equation}
Substituting $Z = \mathcal{C} t$ and $Y = t e^{y-t+q\int_0^t d\tau\,s(\tau)}$ into \eqref{Lagrange} and comparing with the definition \eqref{C:def} of the generating function we arrive at a formal solution 
\begin{equation}
\label{ckt:simple}
c_k(t)= \frac{k^{k-2}}{k!}\, t^{k-1}\,\exp\left[-kt+k q\int_0^t s(\tau)\,d\tau\right]
\end{equation}

\begin{figure}[t]
\includegraphics[width=0.44\textwidth]{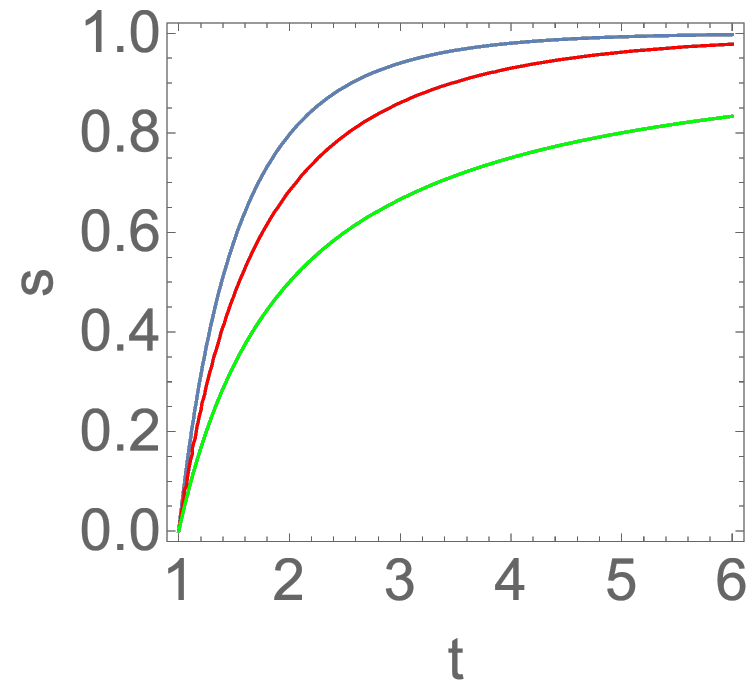}
\caption{The mass $s(t)$ of unicycles determined by Eq.~\eqref{T:sol} in the supercritical  phase, $t>1$. We show $s(t)$ for $p=1, \frac{1}{2}, 0$ (top to bottom). For $p=1$, the mass of unicycles coincides with the mass $g(t)$ of the giant component in the classical random graph process. When $p=0$, the  mass of unicycles is particularly simple: $s=1-1/t$. } 
\label{fig:S}
\end{figure}

To complete the solution \eqref{ckt:simple} we must determine $s(t)$. This can be accomplished by specializing \eqref{y:sol} to $y=0$ and using relation \eqref{def:s}, i.e., $\mathcal{C}(0,t) = 1-s(t)$. We obtain   
\begin{equation}
\label{s:eq}
q\int_1^t d\tau\,s(\tau) = st + \ln(1-s)
\end{equation}
which implicitly determines $s(t)$. Differentiating the integral equation \eqref{s:eq} with respect to time gives an ODE
\begin{equation}
\label{s:ODE}
(1-t+st)\frac{d s}{dt} = p s(1-s)
\end{equation}
Changing $s=s(t)$ to $t=t(s)$ we obtain
\begin{equation}
\label{T:eq}
p\,\frac{dt}{ds} = \frac{1-t}{s} + \frac{1}{1-s}
\end{equation}
Integrating Eq.~\eqref{T:eq} subject to $t(0)=1$ we arrive at the announced result \eqref{T:sol}. 

The quantity $s(t)$ plays a role of an order parameter. The phase transition is continuous [see Fig.~\ref{fig:S}] and mean-field in nature as manifested by the asymptotic  behavior of the order parameter, $s\propto (t-1)$, near the transition. More precisely, 
\begin{subequations}
\begin{equation}
\label{s:small}
s = (1+p)T- \tfrac{(1+p)^3}{1+2p}T^2 +  \tfrac{(1+p)^4(1+4p+2p^2)}{(1+2p)^2(1+3p)}T^3  + \ldots
\end{equation}
when $T=t-1\downarrow 0$. When $t\to\infty$, 
\begin{equation}
\label{s:large}
s = 1-\mathcal{E}-(t-1)\mathcal{E}^2 + \ldots
\end{equation}
\end{subequations}
with $\mathcal{E}(t)$ defined by Eq.~\eqref{E:def}. 

We compute the integral in \eqref{ckt:simple} using the same trick as in the computation of the integral \eqref{integral}. Namely we treat $\sigma=s(\tau)$ as an integration variable and use
\begin{equation}
\label{tau:eq}
p\,\frac{d\tau}{d\sigma} = \frac{1-\tau}{\sigma} + \frac{1}{1-\sigma}
\end{equation}
following from Eq.~\eqref{T:eq}. We get
\begin{eqnarray*}
p\int_0^t \sigma\, d\tau &=& \int_0^s \sigma\left[ \frac{1-\tau}{\sigma} + \frac{1}{1-\sigma}\right] d\sigma  \nonumber\\
&=& \int_0^s \left[\frac{1}{1-\sigma}-\tau \right] d\sigma  \nonumber\\
&=& -\ln(1-s) - t s+\int_0^t \sigma\, d\tau
\end{eqnarray*}
from which we deduce
\begin{equation}
\label{int}
q\int_0^t \sigma\, d\tau = ts+\ln(1-s)
\end{equation}
Using \eqref{int} we recast \eqref{ckt:simple} into
\begin{equation}
\label{ckt:final}
c_k(t)= \frac{k^{k-2}}{k!}\, t^{k-1}\,(1-s)^k\,e^{-kt(1-s)}
\end{equation}

The densities \eqref{prod-densities} in classical random graphs remain smooth throughout the evolution, thereby hiding the phase transition. In SRGs, the densities \eqref{ckt:final} undergo a jump in the second derivative $\frac{d^2 c_k}{dt^2}$ at $t=1$ for all $k\geq 1$ when $p<1$. For instance ($T=t-1$),
\begin{equation}
c_1=e^{-1}\times
\begin{cases}
1-T+(1-p^2/2)T^2+\ldots   & T\downarrow 0\\
1-T+T^2/2+\ldots               & T\uparrow 0
\end{cases}
\end{equation}

\begin{figure}[t]
\includegraphics[width=0.44\textwidth]{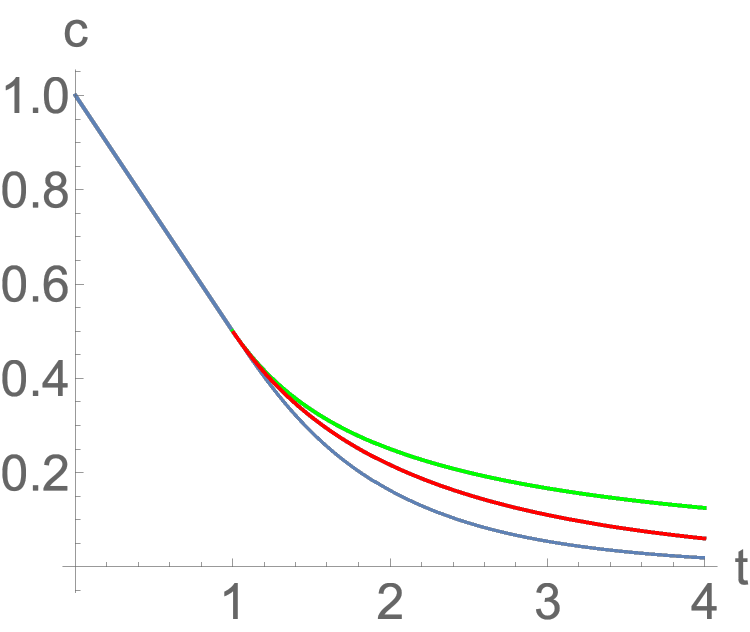}
\caption{Time dependence of the density of trees $c(t)$ determined by Eq.~\eqref{simple-total}. We plot $c(t)$ in the supercritical  phase 
for $p=1, \frac{1}{2}, 0$ (bottom to top).  When $p=0$, the total density of trees is particularly simple: $c=\frac{1}{2t}$. In the subcritical phase, $c=1-t/2$ independently on $p$. } 
\label{fig:C}
\end{figure}

The total density of trees is given by
\begin{equation}
\label{simple-total}
c(t) =
\begin{cases}
1-\frac{t}{2}                                               & t\leq 1\\
1-\tfrac{1}{2}\big(1+s^2\big)t + (t-1)s        & t>1
\end{cases}
\end{equation}
In the subcritical  phase, the total density of trees is the same as in classical random graphs as the governing equations coincide [see Fig.~\ref{fig:C}]. To derive $c(t)$ in the supercritical phase, we sum Eqs.~\eqref{post} and use \eqref{def:s} to find
\begin{equation}
\frac{d c}{d t} =  - \frac{1}{2}\,(1-s)[1+(2p-1)s]
\end{equation}
Integrating this equation we obtain 
\begin{equation}
\label{ct-sol-formal}
c = 1-\frac{t}{2} + q\int_0^t \sigma\, d\tau + \left(p-\frac{1}{2}\right)\int_0^t \sigma^2\, d\tau
\end{equation}
in the supercritical  phase. We already computed the first integral on the right-hand side of \eqref{ct-sol-formal}, see Eq.~\eqref{int}. Employing the same approach, we determine the second integral. Combining these results yields the total density \eqref{simple-total}. In the long time limit, the density decays exponentially, namely as 
\begin{equation}
\label{c:large}
c = \mathcal{E} + \left(\frac{t}{2}-1\right)\mathcal{E}^2  + \ldots
\end{equation}
with $\mathcal{E}$ given by \eqref{E:def}. 

\begin{figure}[t]
\includegraphics[width=0.44\textwidth]{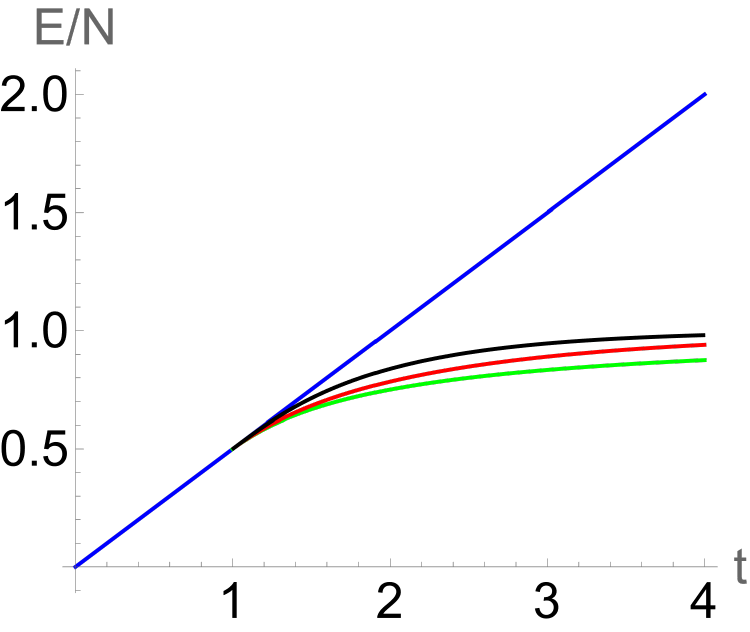}
\caption{Time dependence of the ratio of the total number of edges to the total number of vertices. For classical random graphs, $\frac{E}{N} = \frac{t}{2}$ for all $t>0$, and in the subcritical  phase, the ratio is the same for SRGs independently on the parameter $p$. In the supercritical  phase, the ratio depends on $p$ as illustrated in the figure for $p=1, \frac{1}{2}, 0$ (top to bottom).  When $p=0$, the ratio is particularly simple: $\frac{E}{N}=1-\frac{1}{2t}$. } 
\label{fig:E/N}
\end{figure}

The ratio of the total number of edges to $N$ can be expressed via the density of trees (see Fig.~\ref{fig:E/N})
\begin{equation}
\frac{E}{N} = \sum_{k\geq 1}(k-1)c_k+s = 1-c
\end{equation}

\begin{figure}[t]
\includegraphics[width=0.44\textwidth]{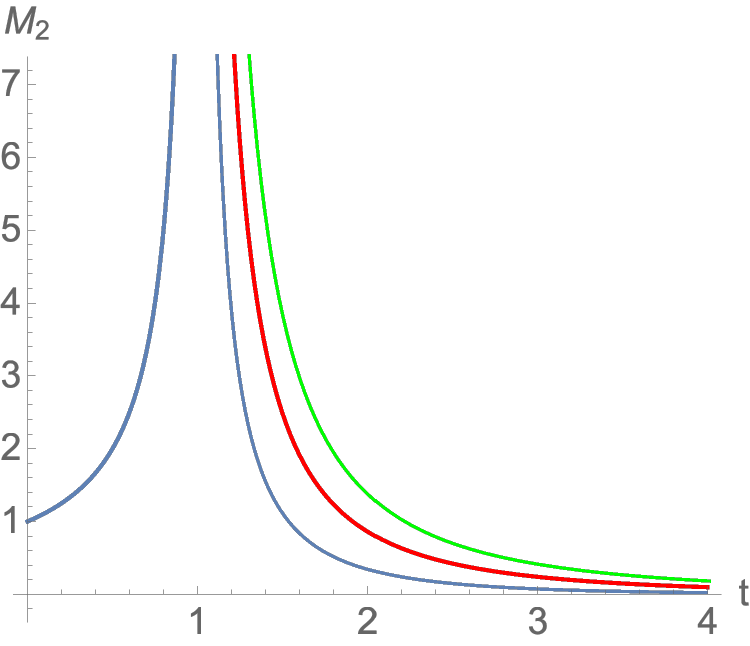}
\caption{Time dependence of the second moment determined by Eq.~\eqref{M2-sol}. We plot the second moment in the supercritical  phase 
for $p=1, \frac{1}{2}, \frac{1}{3}$ (bottom to top). In the subcritical  phase $M_2=1/(1-t)$.  } 
\label{fig:M23}
\end{figure}

To determine the total number of unicycles, we will need the second moment of the distribution of trees. This moment can be expressed via the generating function:
\begin{equation}
\label{M2-GF}
M_2 = \partial_y \mathcal{C}|_{y=0}
\end{equation}
Using \eqref{int} we re-write \eqref{y:sol} as
\begin{equation}
\label{GF:sol}
\mathcal{C}\,e^{-\mathcal{C} t} = (1-s)\, e^{y-t(1-s)}
\end{equation}
We take the logarithm of \eqref{GF:sol}, differentiate with respect to $y$, and set $y=0$ to give
\begin{equation}
\label{M2-sol-gen}
M_2(t) = \left[\frac{1}{\mathcal{C}(0,t)}-t\right]^{-1}
\end{equation}
Recalling that $\mathcal{C}(0,t)=1$ in the subcritical  phase and $\mathcal{C}(0,t)=1-s$ in the supercritical  phase we obtain 
\begin{equation}
\label{M2-sol}
M_2(t) = 
\begin{cases}
(1-t)^{-1}                     & t<1\\
\frac{1-s}{1-t(1-s)}       & t>1
\end{cases}
\end{equation}
Note the similarity with Eq.~\eqref{prod-M2} giving $M_2$ for classical random graphs. When $p=1$, we have $s=g$ and \eqref{M2-sol} reduces to \eqref{prod-M2}. In Fig.~\ref{fig:M23}, we plot the second moment at three values $p=1, \frac{1}{2}, \frac{1}{3}$ of the parameter. 

Using \eqref{s:small}--\eqref{s:large} one finds extremal behaviors
\begin{equation}
\label{M2-post}
M_2(t) = 
\begin{cases}
\frac{1}{p(t-1)}      & t\downarrow 1 \\
\mathcal{E}(t)       & t\uparrow \infty
\end{cases}
\end{equation}
in the supercritical  phase.

\subsection{Unicycles}
\label{subsec:uni}

The number of trees decreases, and at some time $t_\text{jam}$, the last tree disappears, and only unicycles remain. At this moment, the graph freezes. The average number of trees $Nc(t)$ approaches to $N\mathcal{E}(t)$ in the large time limit [cf. \eqref{c:large}]. Combining the criterion $Nc(t_\text{jam})\sim 1$ with \eqref{E:def} we arrive at the asymptotic \eqref{jamming} for the jamming time. The logarithmic scaling with $N$ is certainly correct (when $p>0$), while the amplitude in Eq.~\eqref{jamming} should be taken with a grain of salt. Indeed, the number of trees is a random variable. Fluctuations are significant when the average number of trees is small. If fluctuations are comparable to the average, the asymptotic \eqref{jamming} holds. Recalling that fluctuations in the number of trees scale $\sqrt{N}$ when $t=O(1)$ suggest $\sqrt{N}$ as an upper bound for fluctuations near the jamming time, and then from $Nc(t_\text{jam})\sim \sqrt{N}$ we get a lower bound, $(2p)^{-1}\ln N$, for the jamming time. These heuristic arguments advocate for the asymptotic bounds
\begin{equation}
\label{jamming-bounds}
(2p)^{-1}\ln N \leq t_\text{jam} \leq p^{-1} \ln N
\end{equation}

We conjecture that the upper bound is the true asymptotic, so we used it in Eq.~\eqref{jamming}. Justifying or disproving Eq.~\eqref{jamming} is left for future. 

Arguments in favor of the scaling law \eqref{U-jam} for the final number of unicycles are also heuristic but stronger than the arguments in favor of  Eq.~\eqref{jamming} since the dominant contribution to $U(t_\text{jam})$ is gathered at times $t=O(1)$ far below the jamming time where fluctuations are asymptotically negligible, so the arguments are more solid and the amplitude in Eq.~\eqref{U-jam} appears exact. The begin by noting that for the SRG process, Eq.~\eqref{U-M2} remains valid in the supercritical  phase: Unicycles cannot disappear, and they are born with rate $\frac{1}{2}M_2$. The exact expressions \eqref{M2-sol} for the second moment are valid for infinite graphs. For finite graphs, we do not know $M_2(t)$ in the scaling window, $|t-1|\sim N^{-1/3}$, but we know the magnitude of the second moment: $M_2\sim N^{1/3}$. To gain insight, let us take a test expression for the second moment
\begin{equation}
\label{M2-test}
M_2^\text{test}(t) = 
\begin{cases}
\frac{1}{\sqrt{(1-t)^2+\epsilon^2}}                   & t<1\\
\frac{1}{\sqrt{p^2(t-1)^2+\epsilon^2}}             & 2>t>1
\end{cases}
\end{equation}
with $\epsilon\sim  N^{-1/3}$. This function is continuous, has the correct width of the scaling window, and provides an excellent approximation in the $t<1$ regime away from the scaling window and in the $\epsilon\ll t-1\ll 1$ regime. The gain of $U$ during the time interval $2<t<t_\text{jam}$ is of the order of unity  as $M_2(t)$ decays exponentially with time [see Eq.~\eqref{M2-post} and Fig.~\ref{fig:M23}]. The choice $t=2$ of the cutoff time is asymptotically irrelevant; any cutoff time exceeding $t=1$ suffices---the dominant contribution comes from a small region containing the scaling window. Thus 
\begin{equation}
\label{test}
\frac{1}{2}\int_0^2 dt\,M_2^\text{test}(t)
\end{equation}
estimates $U_\text{jam}$. Using \eqref{M2-test} and computing the integral in \eqref{test} gives \eqref{U-jam} in the leading order. 

The specific form \eqref{M2-test} plays little role, the important feature is the logarithmic divergence  
\begin{equation*}
\int_0^t dt'\,M_2(t')\simeq \ln\frac{1}{1-t}\,, \quad \int_t^{t_\text{jam}} dt'\,M_2(t')\simeq  \frac{1}{p}\ln\frac{1}{t-1}
\end{equation*}
of the integrals near $t=1$ which follows from \eqref{M2-sol} and \eqref{M2-post}. For finite $N$, we must use $|t-1|\sim N^{-1/3}$ to find the integrals in the $(0,1)$ and $(1,t_\text{jam})$ regions. Thanks to the logarithmic divergence we obtain an asymptotically exact leading behavior \eqref{U-jam}.

As another piece of evidence in favor of the prediction \eqref{U-jam} we notice that it reduces to 
\begin{equation}
\label{U-jam-1/2}
U_\text{jam} \simeq \frac{1}{2}\,\ln N
\end{equation}
for $p=\frac{1}{2}$. An `intermediate' SRG process with $p=\frac{1}{2}$ appears in a recent study of a parking process on Cayley trees \cite{Curien23}. This SRG process has an intriguing connection with random maps. A map of a set to itself, $f: S\to S$, can be represented by a graph with edges $(x, f(x))$ for all $x\in S$. This graph decomposes into  maximal connected components (communities). If maps are uniformly chosen among $N^{-N}$ possible maps ($N=|S|$), the distributions of the number of communities and the number of unicycles in the jammed state of the intermediate SRG process are the same. In particular, the average number of communities in random maps is well-known \cite{RM-Kruskal,RM-Harris,Flajolet90}, and it grows with $N$ according to \eqref{U-jam-1/2}. 

Other results about communities in random maps can be restated in terms of unicycles in the SRG process with $p=\frac{1}{2}$. Restating the prediction for the probability to have a single community \cite{RM-Katz} we obtain the probability to end up with one unicycle 
\begin{subequations}
\begin{equation}
\label{U=1}
\text{Prob}[U_\text{jam}=1] = \frac{(N-1)!}{N^N}\sum_{n=0}^{N-1} \frac{N^n}{n!}\simeq \sqrt{\frac{\pi}{2N}}
\end{equation}
The probability to have the maximal number of communities is
\begin{equation}
\label{U=N}
\text{Prob}[U_\text{jam}=N] = \frac{1}{N^N}
\end{equation}
\end{subequations}

Let $K$ be the size of the largest unicycle at jamming. The distribution $\Psi_p(K,N)$ acquires the scaling form
\begin{equation}
\label{KN-scaling}
\Psi_p(K,N) = N^{-1}\psi_p(\kappa)
\end{equation}
when $K\to\infty$ and $N\to\infty$ with $\kappa = K/N$ kept finite. The distribution $\psi_p(\kappa)$ has singularities at $\kappa=1/m$ with $m=2,3,\ldots$ weakening as $m$ increases. The origin of singularities is easy to appreciate: One unicycle may have size exceeding $N/2$, i.e., in the $\frac{1}{2}<\kappa<1$ range; two unicycles may have size in the $\frac{1}{3}<\kappa<\frac{1}{2}$ range; etc. This infinite set of singularities prevents an analytical determination of $\psi_p(\kappa)$. For the intermediate SRG process, Eqs.~\eqref{U=1}--\eqref{U=N} imply the asymptotic behaviors
\begin{equation}
\label{psi-asymp}
\psi_\frac{1}{2}(\kappa) \propto
\begin{cases}
(1-\kappa)^{-1/2}    & \kappa\uparrow 1\\
\kappa^{1/\kappa}    & \kappa\downarrow 0
\end{cases}
\end{equation}

The scaling form \eqref{KN-scaling} implies a linear scaling with $N$ of the average size of the largest unicycle 
\begin{equation}
\langle K\rangle = \lambda(p) N, \quad \lambda(p) = \int_0^1 d\kappa\,\kappa \psi_p(\kappa)
\end{equation}

For classical random maps, the size distribution of the largest community was probed numerically in \cite{RM-Derrida}, and the average size of the largest community is analytically known \cite{Flajolet90} giving $\lambda(\frac{1}{2})=0.757\,823\ldots$.

Consider now  the average number of unicycles $U_k(t)$. For SRG processes with arbitrary $p>0$, it is difficult to probe $U_k(t)$ analytically even in the subcritical phase where it satisfies 
\begin{equation}
\label{Ukt-simple}
\frac{dU_k}{dt}=\frac{1}{2}\,k^2c_k+p\sum_{i+j=k}iU_i\,jc_j-p\,k\,U_k
\end{equation}
These equations differ from Eqs.~\eqref{Ukt} only by the factor $p$ in the last two terms on the right-hand side. Unfortunately, Eqs.~\eqref{Ukt-simple} do not admit a simple general solution like the solution \eqref{Ukt-sol} of Eqs.~\eqref{Ukt} in the case of classical random graphs. 

Equations \eqref{Ukt-simple} are recursive, so one can solve them one by one. The average number of smallest unicycles is 
\begin{equation}
\label{U1p:sub}
U_1 = \frac{e^{-pt}-e^{-t}}{2(1-p)}
\end{equation}
in the subcritical phase. For larger unicycles, the formulas quickly become cumbersome:
\begin{equation*}
U_2 = \frac{(2+p)e^{-2pt}-2pe^{-(1+p)t} - (2-p+4q t)e^{-2t}}{8(1-p)^2} 
\end{equation*}
etc. There is no simple ansatz like \eqref{Ukt-ansatz} fixing the temporal behavior of $U_k(t)$.  

To shed light on the distribution $U_k(t)$ let us look at the behavior of the moments. We have examined the behavior the zeroth moment, $U(t)=\sum_{k\geq 1} U_k(t)$, and know that in the subcritical  phase $U(t)=\frac{1}{2}\,\ln \frac{1}{1-t}$. The first moment, i.e., the average mass of unicycles
\begin{equation}
S(t) = \sum_{k\geq 1} k U_k(t)
\end{equation}
satisfies 
\begin{equation}
\label{S:eq-long}
\frac{dS}{dt} = \sum_{k\geq 1} k \left[ \frac{1}{2}\,k^2c_k+p\sum_{i+j=k}iU_i\,jc_j- p\,k\,U_k \right]
\end{equation}
in the subcritical  phase. Massaging the sums in \eqref{S:eq-long} we arrive at a neat formula
\begin{equation}
\label{S:eq}
\frac{dS}{dt} = \frac{1}{2}\,M_3 + pM_2 S
\end{equation}
We already know $M_2 = \sum_{k\geq 1} k^2 c_k=(1-t)^{-1}$. Similarly one finds $M_3 = \sum_{k\geq 1} k^3 c_k=(1-t)^{-3}$ in the subcritical  phase. Combining these results with \eqref{S:eq} gives
\begin{equation}
\label{S-sol}
S = \frac{1}{2(2-p)}\left[\frac{1}{(1-t)^2} - \frac{1}{(1-t)^{p}}\right]
\end{equation}
applicable in the subcritical  phase. Combining \eqref{S-sol} with $1-t\sim N^{-1/3}$ describing the scaling window we deduce the scaling behavior $S(1)\sim N^{2/3}$ of the average mass of unicycles in the critical regime. In the supercritical phase, $S(t) = N s(t)$ by definition of $s(t)$. 

Higher moments can be probed in a similar way. The exact expressions are cumbersome  already for the second moment, so we just mention the leading asymptotic in the $t\uparrow 1$ limit:
\begin{equation}
\label{n-moment}
\sum_{k\geq 1} k^n U_k(t)\simeq \frac{A_n}{(1-t)^{2n}}
\end{equation}
with amplitudes
\begin{equation}
\label{A12}
A_1= \frac{1}{2(2-p)}\,, \quad A_2= \frac{p}{[2(2-p)]^2}+ \frac{3}{4(2-p)}
\end{equation}
etc. From \eqref{n-moment} we conclude that in the scaling window, the $n^\text{th}$ moment diverges as $N^{2n/3}$. 

\section{Extreme Simple Random Graphs}
\label{sec:extreme}

Many of the above formulas [e.g., Eqs.~\eqref{jamming}--\eqref{U-jam}] become singular when $p=0$. The model with $p=0$ in which unicycles are frozen exhibits different behaviors than the models with $p>0$. Peculiar behaviors occur at a few other values of $p$, e.g., amplitudes \eqref{A12} become singular when $p=2$ indicating that the moments $\sum_{k\geq 1} k^n U_k(t)$ diverge faster than $(1-t)^{-2n}$ as $t\uparrow 1$. However, a natural interpretation of the parameter $p$ as the probability suggests to consider the $p\in [0,1]$ range. We now outline some behaviors for the SRG processes with extreme values of gluing probability. 

\subsection{$p=1$}
\label{subsec:extreme-1}

When $p=1$, Eqs.~\eqref{Ukt-simple} coincide with Eqs.~\eqref{Ukt}. (A singularity in the solution \eqref{U1p:sub} at $p=1$ disappears if one carefully takes the $p\to 1$ limit.) Thus in the subcritical  phase, the solution is given by  \eqref{Ukt-sol}. 

In the supercritical phase
\begin{equation}
\label{Ukt-1}
\frac{dU_k}{dt}=\frac{1}{2}\,k^2c_k+\sum_{i+j=k}iU_i\,jc_j-k\,U_k(1-g)
\end{equation}
where we have taken into account that $s=g$ when $p=1$. 

Consider the smallest unicycles. Solving 
\begin{equation}
\label{U11-eq}
\frac{dU_1}{dt}+U_1(1-g)=\frac{1}{2}\,c_1 = \frac{1}{2}\,e^{-t}
\end{equation}
yields
\begin{equation}
\label{U11-long}
2U_1 e^{t-G(t)} = 1 + \int_1^t dt'\,e^{-G(t')}
\end{equation}
where we shortly write
\begin{equation}
G(t) = \int_1^t d\tau\,g(\tau)
\end{equation}
To compute this integral we treat again $G$ as a function of $g$ rather than $t$. Repeating the steps used in the computation of the integral \eqref{integral} we obtain
\begin{equation}
G = -\ln(1-g)-\text{Li}_2(g)
\end{equation}
where $\text{Li}_2(g)=\sum_{n\geq 1} \frac{g^n}{n^2}$ is the dilogarithm. 

Similarly we simplify the integral on the right-hand side of \eqref{U11-long} and arrive at the solution
\begin{equation}
\label{U11}
2U_1 e^{t(1-g)+\text{Li}_2(g)} = 1 + \int_0^g dh\,H(h)\,e^{\text{Li}_2(h)}
\end{equation}
in the supercritical phase with
\begin{equation*}
\label{H:def}
H(h) = \frac{1}{h}+\frac{(1-h)\ln(1-h)}{h^2}
\end{equation*}

\begin{figure}[t]
\includegraphics[width=0.44\textwidth]{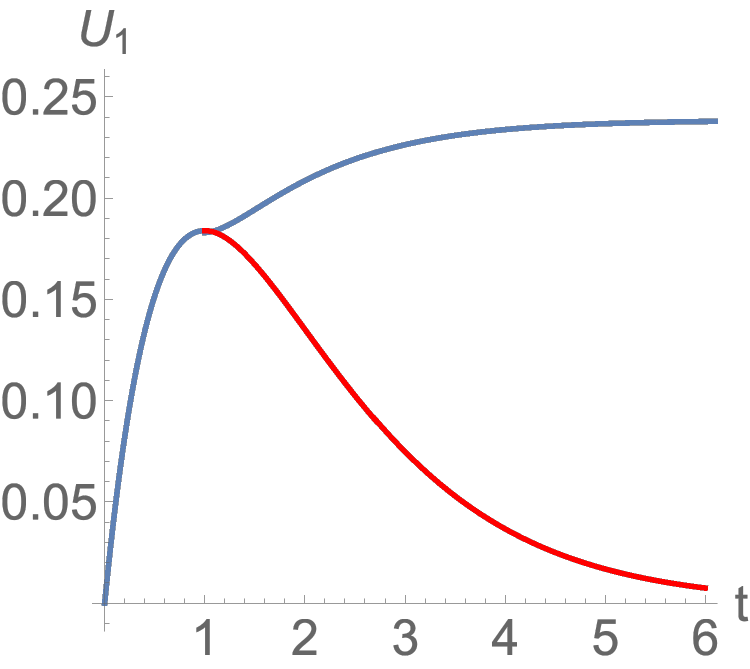}
\caption{Time dependence of the average number of smallest unicycles $U_1(t)$ for SRGs with $p=1$. In the subcritical phase $U_1=\frac{1}{2}te^{-t}$, the same as for the classical random graphs. In the supercritical  phase, $U_1(t)$ is determined by \eqref{U11}. For comparison, we also show  (bottom curve) $U_1(t)$ for classical random graphs. } 
\label{fig:U1}
\end{figure}

Figure \ref{fig:U1} shows that the average number of smallest unicycles is an increasing function of time over the entire evolution, $0<t<t_\text{jam}$. In particular, the final average number of smallest unicycles is   
\begin{equation}
\label{U1-jam}
U_1(t_\text{jam})= \frac{1 + \int_0^1 dh\,H(h)\,e^{\text{Li}_2(h)}}{2e^{\pi^2/6}}\approx 0.23898433
\end{equation}
Thus, at least one smallest unicycle arises with a finite probability. In contrast, the average number of smallest unicycles at the same time $t_\text{jam}\simeq \ln N$ is inverse of the number of vertices, $U_1(t_\text{jam})\sim N^{-1}$, for classical random graphs.

\subsection{$p=0$}
\label{subsec:extreme-0} 

The densities in the subcritical phase are given by \eqref{prod-densities} independently on $p$. The densities in the supercritical phase are given by \eqref{ckt:final} with $s=1-1/t$, the latter result follows from \eqref{T:sol} in the $p\to 0$ limit. Thus, in the supercritical phase ($t>1$), we arrive at
\begin{equation}
\label{S-densities}
c_k(t) = \frac{k^{k-2}e^{-k}}{k!}\, t^{-1}       
\end{equation}
known as the Stockmayer solution \cite{Stockmayer43}, see \cite{Ziff80,Leyvraz-rev}. The total cluster density is particularly simple
\begin{equation}
\label{Stock-total}
c(t) =
\begin{cases}
1-\frac{t}{2}           & t\leq 1\\
\frac{1}{2t}            & t>1
\end{cases}
\end{equation}

The governing equations
\begin{equation}
\label{Ukt-0}
\frac{dU_k}{dt}=\frac{1}{2}\,k^2c_k
\end{equation}
for the average number of unicycles are valid throughout the entire evolution. In the subcritical phase, we insert the densities of trees \eqref{prod-densities} into \eqref{Ukt-0} and integrate to give 
\begin{equation}
\label{Ukt-0-sol}
U_k(t)=\frac{1}{2k}\,e^{-kt}\sum_{n\geq k}\frac{(kt)^{n}}{n!}
\end{equation}
This solution is different from the solution \eqref{Ukt-sol} for classical random graphs. At the critical point 
\begin{equation}
\label{Uk-0-sol}
U_k(1)=\frac{1}{2k}\,e^{-k}\sum_{n\geq k}\frac{k^{n}}{n!}
\end{equation}
differs from $U_k(1)$ for classical random graphs, Eq.~\eqref{Uk-gel}, but the large $k$ behaviors are the same: $U_k(1)\simeq \frac{1}{4k}$ when $k\gg 1$. 

In the supercritical phase, the densities of trees are given by \eqref{S-densities} which we insert into \eqref{Ukt-0} and deduce
\begin{equation}
\label{Ukt-S}
U_k(t)=U_k(1)+\frac{(k/e)^k}{2k!}\,\ln t
\end{equation}
with $U_k(1)$ given by \eqref{Uk-0-sol}. 

We now outline chief properties of the jammed state. Using $c=(2t)^{-1}$ together with criterion $Nc(t_\text{jam})\sim 1$ we arrive at a linear scaling 
\begin{equation}
\label{jamming-0}
t_\text{jam} \sim N
\end{equation}
of the jamming time. (For the models with $p>0$, the scaling is logarithmic.) Using Eqs.~\eqref{Ukt-S}--\eqref{jamming-0} and the Strirling formula we deduce the large $k$ asymptotic 
\begin{equation}
\label{Ukt-jam}
U_k(t_\text{jam})\simeq (8\pi k)^{-\frac{1}{2}}\,\ln N
\end{equation}
at the jamming time. There are only unicycles at the jamming time. Therefore
\begin{equation}
\label{NK}
N = \sum_{k=1}^N k U_k(t_\text{jam}) \sim \ln N \sum_{k\lesssim K} k^\frac{1}{2} \sim  K^\frac{3}{2} \ln N 
\end{equation}
implying that Eq.~\eqref{Ukt-jam} is applicable when $1\ll k \lesssim K$ with
\begin{equation}
\label{K}
K \sim \left(\frac{N}{\ln N}\right)^\frac{2}{3}
\end{equation}
Similarly to the estimating of the sum in \eqref{NK} we deduce 
\begin{equation}
U_\text{jam} = \sum_{k=1}^N U_k(t_\text{jam}) \sim \sum_{k\lesssim K} k^{-\frac{1}{2}} \ln t_\text{jam}  \sim  K^\frac{1}{2} \ln N
\end{equation}
which we combine with \eqref{K} to obtain
\begin{equation}
\label{U-jam-0}
U_\text{jam} \sim N^\frac{1}{3}\,(\ln N)^\frac{2}{3}
\end{equation}

Thus, in the model with frozen unicycles ($p=0$), the average number of unicycles in the jammed state exhibits a peculiar scaling \eqref{U-jam-0}. The size $K$ is smaller by $(\ln N)^{2/3}$ than the maximal size of unicycles in the jammed state: The components of size $N^{2/3}$ at the phase transition point are unicycles, or will become unicycles.

\section{Fluctuations}
\label{sec:fluct}

Fluctuations are relatively small  in large systems, and they are traditionally investigated in the realm of the van Kampen expansion \cite{VanKampen}. Van Kampen expansions have been used in the analyses of various reaction processes \cite{VanKampen,KRB} including aggregation \cite{Lushnikov78,Spouge85b,Ernst87a,Ernst87b}. In our problem, the total number of trees is expected to be the sum of the linear in $N$ deterministic contribution and proportional to $\sqrt{N}$ stochastic contribution:
\begin{align}
\label{trees}
T(t) = N c(t) + \sqrt{N} \xi(t)
\end{align}
The criterion $Nc(t_\text{jam})\sim 1$ for estimating the jamming time is rather naive. The criterion 
\begin{equation}
\label{crit}
Nc(t_\text{jam})\sim \sqrt{N} v(t_\text{jam}), \quad  v(t)=\sqrt{\langle\xi^2(t)\rangle}
\end{equation}
is better, albeit still non-rigorous as it assumes  that the $\sqrt{N}$ scaling of fluctuations holds till the very end.

For the SRG processes with $p>0$, the density of trees decays exponentially, and the criterion \eqref{crit} leads to the logarithmic scaling, and only the amplitude can differ. For instance, if $v(t_\text{jam})=O(1)$, the amplitude is twice smaller than predicted by naive criterion. As we argued in Sec.~\ref{subsec:uni}, these two estimates provides asymptotic bounds \eqref{jamming-bounds}. 

The naive criterion is difficult to justify, but it often leads to asymptotically exact results. For classical random graphs, the naive criterion $Nc(t_\text{cond})\sim 1$, gives the leading behavior with correct amplitude, $t_\text{cond}\simeq \ln N$, of the condensation time when the graph becomes connected \cite{Bollobas,Hofstad,Frieze}. 

For the model with frozen unicycles, $p=0$, fluctuations are potentially more important than for models with $p>0$ where they only affect the amplitude \eqref{jamming}. Indeed, the linear scaling \eqref{jamming-0} of the jamming time is based on the naive criterion. If $v(t_\text{jam})=O(1)$, the criterion \eqref{crit} leads to diffusive, $t_\text{jam} \sim \sqrt{N}$, rather than linear scaling of the jamming time. 

The jammed states in the model with frozen unicycles are remarkably similar to jammed states \cite{Wendy21,Wendy23,Sergey23} discovered in a few addition-fragmentation processes. These jammed states are known as super-cluster states, as clusters tend to be large and the total number of clusters is non-extensive (scales sub-linearly with $N$). 

In the model with frozen unicycles, we have not computed the variance $v(t)$, so the precise scaling of the jamming time remains unknown. The formula \eqref{Ukt-S} remains valid before the jamming, and \eqref{Ukt-jam} with $\ln t_\text{jam}$ instead of $\ln N$. We anticipate $t_\text{jam}\sim N^b$, which would modify \eqref{Ukt-jam} by a numerical factor. This factor is irrelevant in deriving \eqref{U-jam-0}--\eqref{K}. The scaling laws \eqref{U-jam-0}--\eqref{K} suggest a similar scaling law 
\begin{equation}
\label{jamming-log}
t_\text{jam} \sim N^b (\ln N)^\beta
\end{equation}
for the jamming time.

\section{Concluding Remarks}
\label{sec:remarks}

Classical random graphs undergo a percolation transition at $t_c=1$ when the giant component is born. The mass of the giant component, i.e., the fraction of vertices belonging to the giant component, plays the role of the order parameter in the supercritical phase $t>t_c=1$. Classical random graphs subsequently condense into a single component, i.e., the giant component engulfs all finite components. This condensation transition occurs at $t_\text{cond}\simeq \ln N$. The densification of classical random graphs continues forever or until the graph becomes complete if loops and multiple edges are forbidden. 

Simple random graphs (SRGs) evolve similarly to classical random graphs with the constraint that the formation of complex components is forbidden. (A complex component has a negative Euler characteristic.)  The SRGs undergo a phase transition at $t_c=1$. When $t\leq 1$, the evolution of SRGs is asymptotically identical to classical random graphs where only a few complex components may arise. In the supercritical phase, $t>1$, the mass of unicycles is an order parameter behaving similarly to the mass of the giant component for classical random graphs. A crucial difference between classical random graphs and SRGs is that for $t>1$, the former contain a single giant component, while SRGs contain many macroscopic unicycles (i.e., unicycles with size proportional to $N$). Moreover, the size of the giant component is a self-averaging random quantity---fluctuations around the average $gN$ scale as $\sqrt{N}$. For the SRGs, the size of the largest unicycle at jamming is a non-self-averaging random quantity with non-trivial distribution \eqref{KN-scaling}. 

We have studied a class of SRG processes depending on parameter $p\in [0,1]$. Some properties of random graphs emerging in models with $p=0$ and $p=1$ resemble classical results about aggregation processes with the merging rate of components proportional to the product of their masses. The interpretation of mathematically identical formulas is different, however. For instance, the mass of unicycles for the SRG process with frozen unicycles ($p=0$) coincides with the gel mass in the Stockmayer model of gel formation \cite{Stockmayer43}. The mass of unicycles for the SRG process with $p=1$ is the same as the gel mass in the Flory model of gelation \cite{Flory41, Flory53}, or equivalently, the mass of the giant component in classical random graphs. An `intermediate' SRG process with $p=\frac{1}{2}$ appears in a recent study of a parking process on Cayley trees \cite{Curien23}. The behavior of this SRG process, particularly in the scaling window around the phase transition point, has been analyzed in Ref.~\cite{Curien23}. 

Jamming is a distinctive feature of the SRG processes. The evolution of a simple random graph freezes when trees disappear, and only unicycles remain. An analytical determination of the distribution of the unicycles at the jamming time for models with $p>0$ is an open issue. We argued that the jamming time and the average number of unicycles in a jammed state scale logarithmically with system size, Eqs.~\eqref{jamming}--\eqref{U-jam}. The arguments leading to these scaling laws are heuristic since we ignored fluctuations. The logarithmic scaling laws seem correct. The amplitude in \eqref{jamming} is a bit questionable and could be just an upper bound, see \eqref{jamming-bounds}. Arguments in favor of the amplitude in \eqref{U-jam} are more solid. 

The SRG process with frozen unicycles ends up in intriguing jammed states resembling super-cluster jammed states \cite{Wendy21, Wendy23, Sergey23} found in aggregation-fragmentation processes. To compute the jamming time, one should  know the variance $\langle\xi^2(t)\rangle$.  We argued, however, that the scaling laws \eqref{U-jam-0}--\eqref{K} derived without knowing the jamming time are asymptotically exact. Relying on the form of the scaling laws \eqref{U-jam-0}--\eqref{K}, we conjectured the scaling law \eqref{jamming-log} for the jamming time. 

Apart from non-extensive super-cluster jammed states \cite{Wendy21,Wendy23,Sergey23}, non-extensive steady states also appear in some aggregation-fragmentation processes \cite{BK08}. The emergence of such non-thermodynamic behaviors often requires tuning the parameters, like setting $p=0$ in the SRG process. Still, non-thermodynamic behaviors appear widespread. The theoretical analysis is challenging since fluctuations often play a decisive role, and computing fluctuations in a system with numerous interacting cluster species tends to be impossible. 

Simple and classical random graphs undergo a `mean-field' continuous phase transition. More exotic phase transitions can also arise in random graph processes. It would be interesting to modify these processes, forbidding the creation of complex components or imposing other constraints. One can then examine the robustness of the phase transition. For instance, many growing networks \cite{BKT-Strogatz,KR01-Rodgers,BKT-Sergey,KR02-protein,KR02-BKT,BKT-Bauer-a,BKT-Bauer-b,BKT-Derrida,BK07} exhibit an infinite-order percolation transition of the Berezinskii-Kosterlitz-Thouless type. Unicycles in such networks are also worth investigating. 

\bigskip\noindent
I am grateful to N. Curien and M. Tamm for comments.

\bibliography{references-graphs}

\end{document}